\documentclass[11pt]{article}

\usepackage{a4wide}

\usepackage{amssymb}
\usepackage{amsthm}
\usepackage{enumerate}
\usepackage{amsmath}

\usepackage{harvard}

\newtheorem{thm}{Theorem}
\newtheorem*{defi}{Definition}
\newtheorem{cor}{Corollary}
\newtheorem{lem}{Lemma}
\newtheorem{proposition}{Proposition}
\newtheorem{rem}{Remark}
\newtheorem*{thm*}{Theorem}

\newcommand{\N}{\mathbb N}
\renewcommand{\L}{\mathcal L}
\newcommand{\Z}{\mathbb Z}
\newcommand{\R}{\mathbb R}
\newcommand{\T}{\mathbb T}
\newcommand{\dif}{{\rm d}}
\newcommand{\ind}{1}

\newcommand{\F}{{\mathcal F}}
\renewcommand{\H}{{\mathcal H}}
\newcommand{\A}{{\mathcal A}}
\newcommand{\B}{{\mathcal B}}

\newcommand{\esp}{\mathbb E}

\newcommand{\cqfd}{\hfill $\Box$}
\newcommand{\Cov}{{\rm Cov}}
\newcommand{\Var}{{\rm Var}}
\newcommand{\cla}{\stackrel{{\cal D}}{\longrightarrow}}
\newcommand{\Cl}{{\mathcal C}^1}

\citationstyle{agsm}
\begin{document}

 \title{\sc{Empirical Processes of Multidimensional Systems with Multiple Mixing Properties}}
 \author{
Herold Dehling\thanks{Fakult\"at f\"ur Mathematik, 
Ruhr-Universit\"at Bochum, Universit\"atsstra\ss e 150, 44780 Bochum, Germany;
e-mail: herold.dehling@ruhr-uni-bochum.de; Research supported by the German Science Foundation, 
Sonderforschungsbereich 823}
\and
Olivier Durieu\thanks{Laboratoire de Math\'ematiques et Physique Th\'eorique, UMR 6083 CNRS, Universit\'e Fran\c{c}ois Rabelais de Tours, Parc de Grandmont, 37200 Tours, France;
e-mail: olivier.durieu@lmpt.univ-tours.fr}
}

\maketitle

\begin{abstract}
We establish a multivariate empirical process central limit theorem for 
stationary $\R^d$-valued stochastic processes $(X_i)_{i\geq 1}$ under very 
weak  conditions concerning the dependence structure of the process. As an
application we can prove the empirical process CLT for ergodic torus 
automorphisms. Our results also apply to Markov chains and dynamical systems
having a spectral gap on some Banach space of functions. Our proof uses a 
multivariate extension of the techniques introduced by \citeasnoun{DehDurVol09}
 in the univariate case. As an important technical ingredient,
we prove a $2p$-th moment bound for partial sums in multiple mixing systems.

\medskip


\medskip

\end{abstract}

\section{Introduction and Statement of Main Results}
Let $(X_i)_{i\geq 1}$ be an $\R^d$-valued stationary stochastic process with
multivariate marginal distribution function $F(t)=P(X_1\leq t)$, $t\in \R^d$.
We define the empirical distribution function and the empirical process by
\begin{eqnarray*}
 F_n(t)&:=& \frac{1}{n} \#\{1\leq i\leq n: X_i\leq t  \} \\
 U_n(t)&:=& \sqrt{n}(F_n(t)-F(t)),
\end{eqnarray*}
$t\in \R^d$. Here "$\leq$" denotes the coordinate-wise ordering, i.e.
$(t_1,\ldots,t_d)\leq (s_1,\ldots,s_d)$ if and only if $t_i\leq s_i$ for all
$i\in \{1,\ldots,d\}$. In this paper, we study weak convergence of the 
empirical process towards a Gaussian process in the space 
$D([-\infty,\infty]^d)$. We make very weak assumptions concerning the 
dependence structure of the underlying process $(X_i)_{i\geq 1}$. 
Effectively, we require a multiple mixing condition, the central limit theorem
for partial sums of a restricted class of functions and a condition on the modulus
of continuity of the distribution function $F$.
Our results apply to Markov chains and dynamical systems whose transfer 
operator has a spectral gap on some Banach space of functions. 
As most significant application, 
we can establish the multivariate empirical process CLT for
non-hyperbolic ergodic torus automorphisms, a system that does not have a 
spectral gap on common spaces of functions.

The study of empirical processes was initiated by Donsker's empirical
process invariance principle \citeaffixed{Don52}{see}, which covered the case of i.i.d. 
$\R$-valued observations.
Donsker's original theorem has been generalized to dependent variables by a
number of authors, starting with work by Billingsley \citeyear{Bil68} who could
establish the empirical process invariance principle for functionals of 
uniformly mixing processes. Billingsley applied this result in his 
investigations of statistical properties of the continued fraction expansion.
Berkes and Philipp \citeyear{BerPhi77} were able to treat the empirical process of
strongly mixing sequences. Borovkova, Burton and Dehling \citeyear{BorBurDeh01} could treat
functionals of absolutely regular processes. Dehling and Taqqu \citeyear{DehTaq89} proved 
an empirical process invariance principle for long-range dependent data.
Dedecker and Prieur \citeyear{DedPri07} proved empirical process invariance principles for processes that satisfy one of the weak dependence conditions introduced earlier by the same authors, see Dedecker and Prieur \citeyear{DedPri05}. These are generalizations of classical mixing coefficients, in a different way than Doukhan and Louhichi \citeyear{DouLou99}; see also the recent
book by Dedecker et al \citeyear{Dedal07}. Wu \citeyear{Wu08} studied empirical processes
in the case that the underlying process can be represented as a
functional of an i.i.d. process $(\epsilon_i)_{i\in \Z}$, 
i.e. $X_i=f((\epsilon_{i-k})_{k\geq 0})$.
 Wu and Shao \citeyear{WuSha04} investigated the empirical process for certain classes of
Markov chains.

In recent years, a lot of research has been devoted to the study of statistical
properties of data arising from dynamical systems. Given a measure 
preserving dynamical system $(\Omega,\F,T,P)$, consider the process
$X_n:=T(X_{n-1})$, $n\geq 1$. When $T$ is a uniformly expanding map of the
unit interval, this process can be represented as a functional of 
an absolutely regular process; see Hofbauer and Keller \citeyear{HofKel82}. Denker
and Keller \citeyear{DenKel86} used this representation in their investigations of the
asymptotic behavior of $U$-statistics when the underlying data arise from
a dynamical system. By combining this representation with coupling ideas,
Borovkova, Burton and
Dehling \citeyear{BorBurDeh01} were able to study the empirical process and more generally,
$U$-processes. 

The spectral gap technique is a  very powerful technique that allows to 
study much larger classes of dynamical systems. Let $Q$ denote the 
Perron-Frobenius operator, defined by $\int f\circ T g dP =\int f Qg dP$,
where $f\in L_\infty$, $g\in L_1$. The spectral gap technique
studies the dynamical system via spectral properties of the Perron-Frobenius
operator, viewed as operator on a suitable invariant subspace of $L_1$. 
The spectral gap technique has been very successfully applied to the study 
of central
limit theorems and large deviations properties;  a survey and a large
number of examples can be found in the monograph by
Hennion and Herv\'e \citeyear{HenHer01}. It is possible to treat empirical processes 
within the framework of the spectral gap technique if the Perron-Frobenius
operator has a spectral gap on the space of functions of bounded variation, see Collet, Martinez and Schmitt \citeyear{ColMarSch04}.
In this case one can directly apply the standard proof of the empirical process
invariance principle, i.e. establish finite dimensional convergence and 
tightness. This is essentially due to the fact that the indicator functions 
are functions of bounded variation. The situation is different when the 
spectral gap property can only be established on a smaller space of functions,
such as Lipschitz functions. Recently, Gou\"ezel \citeyear{Gou08} has given a 
one-dimensional example of such a dynamical system.

Dehling, Durieu and Voln\'y \citeyear{DehDurVol09} have developed a new technique
that is particularly useful when handling data from dynamical systems 
whose Perron-Frobenius operator has a spectral  gap in the space of Lipschitz
functions. Instead of trying to deduce inequalities for bounded variation functions
from those on Lipschitz functions, the authors proposed a proof which only involves
Lipschitz functions. 
In principle, this technique can also be applied to other spaces of functions, not just Lipschitz functions.
The technique can also be applied to the study of Markov processes whose 
Markov operator has a spectral gap in the space of Lipschitz functions.
The new technique uses classical chaining ideas, 
but replaces the indicator functions that are commonly used 
by Lipschitz functions. Dehling, Durieu and Voln\'y \citeyear{DehDurVol09} make two assumptions
concerning the process $(X_i)_{i\geq 0}$. For any Lipschitz functions 
$f:\R\rightarrow \R$, the partial sums $\sum_{i=1}^n f(X_i)$ satisfy the 
central limit theorem and a suitable bound on the $4$-th central moments.
Under these two assumptions and a mild additional
assumption on the modulus of continuity 
of the distribution function $F$, Dehling, Durieu and Voln\'y \citeyear{DehDurVol09} could 
establish the empirical process invariance principle. 

In the present paper, we extend the techniques of Dehling, Durieu and 
Voln\'y \citeyear{DehDurVol09} to multidimensional systems satisfying a multiple mixing 
condition. In the multidimensional case, the 4th moment bounds have to be
replaced by bounds on higher order moments. We establish such a bound for 
multiple mixing systems. The multiple mixing condition has been used in the
study of the statistical properties of dynamical systems, e.g. in the
work of Le Borgne \citeyear{Leb99} and Durieu and Jouan \citeyear{DurJou08} on non-hyperbolic 
torus automorphisms.

\begin{defi}
Let $(\B,\|\,\|)$ be a Banach space of measurable functions $\varphi: \R^d
\longrightarrow \R$. We say that the process $(X_n)_{n\geq 0}$ has a multiple
mixing property with respect to $\B$ if there exist constants $0<\theta<1$
and $r\geq 1$ such that for any $p\in \N\setminus \{0\}$, 
there exist a positive constant $C$ and an integer $\ell$ such that the following 
assertions hold:
for any $i_1,\ldots,i_p\in \N$ and $q\in \{1,\ldots, p\}$, for any $\varphi$ in $\B$ such that $\esp_\nu\left(\varphi\right)=0$ and $\|\varphi\|_\infty\le 1$,
\begin{equation}\label{mm}
|\Cov(\varphi(X_0)\varphi(X_{i_1^*})\dots\varphi(X_{i_{q-1}^*}),
\varphi(X_{i_q^*})\dots\varphi(X_{i_p^*}))|\le C\|\varphi(X_0)\|_r\|\varphi\|^{\ell}P(i_1,\dots,i_p)\theta^{i_q},
\end{equation}
where $P$ is a polynomial function of $p$ variables which does not depend on $\varphi$ and where we use the following notation:
if $(a_n)_{n\ge 1}$ is a sequence of real number, $a_n^*$ is the sum $\sum_{i=1}^n a_i$.
\end{defi}

In this paper, we will work mostly with the Banach space $\H^\alpha$ of 
bounded $\alpha$-H\"older continuous functions, $0<\alpha \leq 1$.  
We define the norm
\begin{equation}\label{Hol}
\|g\|=\|g\|_\infty+\sup_{x\neq y\in \R^d}\frac{|g(x)-g(y)|}{|x-y|^\alpha}.
\end{equation}

\begin{thm}
\label{th:mmepclt}
Let $(X_i)_{i \geq 1}$ be an $\R^d$-valued stationary stochastic process 
satisfying the multiple mixing property with respect to the Banach space
$\H^\alpha$. Assume that for all $\phi \in \H^\alpha$, the partial 
sums $\sum_{i=1}^n \phi(X_i)$ satisfy the central limit theorem, i.e. that
\[
 \frac{1}{\sqrt{n}} \sum_{i=1}^n (\phi(X_i)-E\phi(X_1)) 
 \rightarrow N(0,\sigma^2)
\]
where $\sigma^2=\Var(\phi(X_1))+2\sum_{i=2}^\infty \Cov(\phi(X_i),\phi(X_1))$.
If the modulus of continuity $\omega$ of the distribution function $F$ satisfies the 
condition 
\begin{equation}\label{modulus}
 \omega(\delta)=O(|\log(\delta)|^{-\gamma})\mbox{ for some } \gamma>r,
\end{equation}
where $r$ is given by (\ref{mm}), then the empirical process central limit theorem holds, i.e.
\[ 
(U_n(t))_{ t \in [-\infty,\infty]^d } \stackrel{\mathcal{D}}{\longrightarrow}
(W(t))_{t \in [-\infty,\infty]^d},
\]
where ``$\stackrel{\mathcal{D}}{\longrightarrow}$'' denotes the weak convergence in the Skorohod space $D([-\infty,\infty]^d)$. \\
Here $(W(t))_{t\in[-\infty,\infty]^d}$ is a mean zero Gaussian process with covariance structure
\begin{eqnarray*} 
 E  W (s) \cdot W (t) 
  &=& \Cov ( 1_{(-\infty, s]}(X_0) , 1_{(-\infty, t]} (X_0))\\
  &&+ \sum^\infty_{k=1} \Cov(1_{(-\infty, s]}(X_0), 1_{(-\infty, t]}(X_k))\\
  &&+ \sum^\infty_{k=1} \Cov (1_{(-\infty, s]}(X_k), 1_{(-\infty, t]}(X_0)).
\end{eqnarray*}
Further, almost surely, $(W(t))_{t\in[-\infty,\infty]^d}$ has continuous sample paths.
\end{thm}
Notice that the assumption on the modulus of continuity of $F$ plays a key role in the proof of the theorem.
It allows us to control the indicator functions of the process from the control on the H\"older observables.
As mentioned by the referee of this paper, there are good reasons to believe that condition (\ref{modulus}) is necessary if one only has assumptions on H\"older functions.

The proof of Theorem~\ref{th:mmepclt} will be given in two parts.
In Section~\ref{sec:ip-1} we will establish a  general empirical process 
CLT under the conditions of Theorem~1, but with the multiple mixing replaced
by a $2p$-th moment bound. In Section~\ref{sec:mb} we will show that 
multiple mixing implies this $2p$-th moment bound.

As an application of Theorem~\ref{th:mmepclt} we can establish the empirical
process invariance principle for ergodic torus automorphism. We consider the
torus $\T^d$, identified to $[0,1]^d$ and equipped with the Lebesgue measure, and define the automorphism $T:\T^d\rightarrow \T^d$
by
\[
 T(x) = M\, x \quad \mod 1.
\]
Here $M$ is a $d\times d$ matrix with integer 
entries and $|{\rm det}M|=1$. We assume that the matrix $M$ has no eigenvalue which is a root of
unity which is equivalent to ergodicity. Such torus automorphisms always have at least one eigenvalue of modulus strictly bigger than 1 and then another one of modulus strictly smaller than 1. Thus a part of the action of the automorphism on the torus has some hyperbolicity. The automorphism is called hyperbolic if it has no eigenvalue of modulus one, and quasi-hyperbolic if it has eigenvalues of modulus one.
For more details on torus automorphisms see Lind \citeyear{Lin82}.
Le Borgne \citeyear{Leb99} established the central limit theorem for 
quasi-hyperbolic torus automorphisms. Durieu and Jouan \citeyear{DurJou08} proved the 
empirical process central limit theorem for certain univariate functionals of 
quasi-hyperbolic torus automorphisms, i.e. they considered the 
empirical distribution of the sequence $(f(T^k\, x))_{k\geq 1}$, for
$f:\T\rightarrow \R$. In this paper we can establish the full empirical 
process invariance principle in the same case.

\begin{thm}\label{thmtorus}
Let $T$ be an ergodic automorphism of the $d$-dimensional torus. Then the 
empirical process
\[
 U_n(t)=\sqrt{n}(F_n(t)-\prod_{i=1}^d t_i), \; t=(t_1,\ldots,t_d)\in [0,1]^d
\]
converges in distribution to a centered Gaussian process 
$(W_{t})_{t\in[0,1]^d}$ which has almost surely continuous sample paths.
\end{thm}

\medskip

\section{An invariance principle for the multivariate empirical process}\label{sec:ip-1}

Let $\left(X_n\right)_{n\ge 0}$ be a stationary process with values in $\R^d$, $d\ge 1$.
For $t=(t_1,\dots,t_d)\in\R^d$ and $s=(s_1,\dots,s_d)\in\R^d$, we use the notations
\[
 t+s=(t_1+s_1,\dots,t_d+s_d),
\]
\[
s\le t \Leftrightarrow s_i\le t_i \mbox{ for all }i=1,\dots,d,
\]
\[
s< t\Leftrightarrow s_i< t_i\mbox{ for all }i=1,\dots,d 
\]
and for $s< t$,
\[
 [s,t]=\prod_{i=1}^d [s_i,t_i].
\]

\medskip

As mentionned in the introduction, to generalize the result of \citeasnoun{DehDurVol09} to multivariate processes, we need
higher moment bounds. In this section, we will assume that the process $\left(X_n\right)_{n\ge 0}$ satisfies moment bounds on a given Banach space.
The technique which is developed here is useful in cases where the Banach space does not contain the indicators functions. Our technique will work if the Banach space $\B$ is enough well adapted to approximate indicator functions. Typically, we could work with space of regular functions as the spaces of Lipschitz continuous functions, H\"older continuous functions, or ${\mathcal C}^k$ functions.
Here, to have a link with the applications, we will work with the space $\H^\alpha$ of bounded $\alpha$-H\"older functions, for some fixed $\alpha\in (0,1]$.
This space is equipped with the norm $\|.\|=\|.\|_{\H^\alpha}$ defined in (\ref{Hol})

\medskip

Then, we make two assumptions concerning the process 
$(X_i)_{i\geq 0},$ 

\begin{enumerate}
\item 
For any function $f\in\H^\alpha$, the CLT holds, i.e. $\sum_{i=0}^\infty\Cov(f(X_0),f(X_i))$ converges and
\begin{equation}
 \frac{1}{\sqrt{n}} \sum_{i=1}^{n}(f(X_i)-Ef(X_i))
  \cla N(0,\sigma^2),
\label{eq:clt}
\end{equation}
where $N(0,\sigma^2)$ denotes a normal law with 
mean zero and variance 
\[
\sigma^2=E(f(X_0)-Ef(X_0))^2+2\sum_{i=1}^\infty\Cov(f(X_0),f(X_i)).
\]

\item 
For any $p\ge 1$, a bound on the $2p$ central moments of partial sums of $(f(X_i))_{i\geq0}$, 
$f\in\H^\alpha$ with $E(f(X_0))=0$ and $\|f\|_\infty\le 1$, of the type
\begin{equation}\label{eq:moment}
\esp_\nu\left[\left(\sum_{i=1}^{n}f(X_i)\right)^{2p}\right]\leq
K\sum_{i=1}^p n^i\|f\left(X_0\right)\|_r^i\log^{2p-i}(\|f\|+a)
 \end{equation}
where K is some universal constant, $a>1$ and $r\ge 1$.
\end{enumerate}
In Section \ref{sec:mb}, we will show how this condition is implied by the multiple mixing property.

\medskip

Recall that
the empirical distribution function $(F_n(t),\,t\in\R^d)$ and the 
empirical process $(U_n(t),\,t\in\R^d)$ are defined by
\begin{eqnarray*}
 F_n(t) &= & \frac{1}{n} \sum_{i=1}^n 1_{(-\infty,t]} (X_i),\; t\in\R^d,
 \\
 U_n(t) &=& \sqrt{n} (F_n(t)-F(t)),\; t\in\R^d,
\end{eqnarray*}
where $F$ is the distribution function of $X_0$ and
the modulus of continuity of a function
$f:\R^d\longrightarrow \R$ by
$$
\omega_f(\delta)=\sup\left\{|f(s)-f(t)|\,:\,s,t\in\R^d,|s-t|<\delta\right\}.
$$

\begin{thm}\label{thm1}
Let $(X_i)_{i\ge 0}$ be an $\R^d$-valued stationary random process such that
the conditions (\ref{eq:clt}) and (\ref{eq:moment}) hold. Assume that $X_0$ has a distribution function $F$ satisfying the following condition,
\begin{equation}\label{eq-modulus}
 \omega_F(\delta)=O(|\log(\delta)|^{-\gamma})\mbox{ for some } \gamma>r,
\end{equation}
where $r$ is given by (\ref{eq:moment}). 
Then
\begin{equation*}
 (U_n(t))_{t\in [-\infty,\infty]^d} \stackrel{\mathcal{D}}
 {\longrightarrow} (W(t))_{t\in [-\infty,\infty]^d},
\end{equation*}
where $W(t)$ is a mean-zero Gaussian process with covariances
\begin{eqnarray*} 
 E  W (s) \cdot W (t) 
  &=& \Cov ( 1_{(-\infty, s]}(X_0) , 1_{(-\infty, t]} (X_0))\\
  &&+ \sum^\infty_{k=1} \Cov(1_{(-\infty, s]}(X_0), 1_{(-\infty, t]}(X_k))\\
  &&+ \sum^\infty_{k=1} \Cov (1_{(-\infty, s]}(X_k), 1_{(-\infty, t]}(X_0)).
\end{eqnarray*}
Further, almost surely, $(W(t))_{t\in\R^d}$ has continuous sample paths.
\end{thm}

\begin{rem}
Here, if (\ref{eq:moment}) holds for $r=1$ as in Dehling et al.\citeyear{DehDurVol09}, then in assumption (\ref{eq-modulus}) we can consider every $\gamma>1$. This is an improvement of the corresponding theorem for dimension 1 of Dehling et al.\citeyear{DehDurVol09}. This is a consequence of the fact that we consider $2p$-th moment inequalities ($p\ge 1$) instead of only a $4$-th moment bound.
\end{rem}

\medskip

In this paper, we work with the Skorohod topology on the function space
$D([-\infty,\infty]^d)$, as introduced by Neuhaus \citeyear{Neu71} and Straf \citeyear{Str72}.
In fact, these authors considered the space $D([0,1]^d)$,  but we can easily
extend their definitions since $[0,1]^d$ and $[-\infty,\infty]^d$ are
homeomorphic. Take any homeomorphism 
$\phi:[-\infty,\infty] \rightarrow [0,1]$, e.g.
\[
 \phi(t)=\frac{1}{2} + \frac{1}{\pi} \arctan (t),
\]
and define $\Phi:[-\infty,\infty]^d \rightarrow [0,1]^d$ by
\[
 \Phi(t_1,\ldots,t_d)=(\phi(t_1),\ldots,\phi(t_d)).
\]
The map $\Phi$ induces a map that associates to any function 
$f\in D([0,1]^d)$ the function $f\circ \Phi:[-\infty,\infty]^d\rightarrow \R$.
We define $D([-\infty,\infty]^d)$ as the image of $D([0,1]^d)$ under this map.
Neuhaus \citeyear{Neu71} and Straf \citeyear{Str72} introduced a metric $d_0$ on 
$D([0,1]^d)$ that generates the Skorohod topology and such that
$(D([0,1]^d),d_0))$ is a complete separable metric space. We can naturally
extend $d_0$ to $D([-\infty,\infty]^d)$ by defining for 
$g_1,g_2\in D([-\infty,\infty]^d)$
\[
 \tilde{d}_0 (g_1,g_2)=d_0(g_1\circ \Phi^{-1}, g_2\circ \Phi^{-1}).
\]
In what follows, we will denote the metric on $D([-\infty,\infty]^d)$
also by $d_0$. Note that $d_0$ is bounded by the supremum distance, i.e.
\[
 d_0(g_1,g_2) \leq \sup_{t\in [-\infty,\infty]^d}|g_1(t)-g_2(t)|,
\]
and that $(D([-\infty,\infty]^d),d_0)$ is a complete separable metric space.
Note that the sample paths of the processes arising  in this paper 
are elements of the function
space $D([-\infty,\infty]^d)$, since their limits as any of the arguments 
approach $\infty$ exist. Alternatively one can see this by observing
that e.g. the empirical distribution function $F_n$ of the process
$(X_i)_{i\geq 1}$ is the image under the map defined above 
of the empirical distribution function of the process $(\Phi(X_i))_{i\geq 1}$.

\medskip

{\em Proof of Theorem~\ref{thm1}.}

To prove Theorem~\ref{thm1}, we shall adapt the technique introduced by Dehling et al. \citeyear{DehDurVol09}.
The idea is to replace the indicator functions $1_{(-\infty,t]}(x)$ by approximations in the space $\H^\alpha$.

\medskip

For each $i=1,\dots,d$, we denote by $F_i$ the marginal distribution functions of $X_0$ corresponding to the i-th coordinate. Note that the $F_i$ also verify condition (\ref{eq-modulus}).

Given a partition of $[0,1]$,
\begin{equation*}
0=r_0< \ldots <r_m=1
\end{equation*}
we define
$$
t_{i,j}=F_i^{-1}(r_j)
$$
where $F_i^{-1}$ is given by
$$
F_i^{-1}(t)=\sup\{s\in \R : F_i(s)\leq t\}.
$$
Thus, by continuity of the $F_i$, we have subdivisions
$$
-\infty\le t_{i,0}<\dots<t_{i,m}=+\infty.
$$
If $j=(j_1,\dots,j_d)\in\{0,\dots\,m\}^d$, we set
\[
 t_j=(t_{1,j_1},\dots,t_{d,j_d}).
\]

We introduce the functions $\varphi_j:\R^d\rightarrow\mathbb{R}$, $j=(j_1,\dots,j_d)\in\{1,\dots\,m\}^d$
defined by
\[
\varphi_j(x)=\left\{
\begin{array}{ll}\displaystyle
\prod_{i=1}^{d}\varphi\left(\frac{x_i-t_{i,j_i-1}}{|t_{i,j_i-1}-t_{i,j_i-2}|}\right)&\mbox{ if } (2,\dots,2)\le j\\
0&\mbox{ otherwise}
\end{array}
\right.
\]
with 
\begin{equation}
 \varphi(x)=1_{(-\infty,-1]}(x)-x1_{(-1,0]}(x)
\label{eq:phi}
\end{equation} 
and where we eventually used the convention that $\frac{1}{\infty}=0$.\\
Note that $\varphi$ is a $\alpha$-H\"older function on $\R$ (for all $\alpha\in (0,1]$).
The function $\varphi_j$ will serve as a $\H^\alpha$-approximation to the 
indicator function $1_{(-\infty,t_{j-1}]}$.

\medskip

Now, we introduce the process
\begin{eqnarray*}
 F_n^{(m)}(t) 
 &=& \frac{1}{n} \sum_{i=1}^{n} 
   \sum_{j\in\{1,\dots\,m\}^d} 1_{[t_{j-1}, t_j)}(t) \varphi_j(X_i)\\
 &=& \sum_{j\in\{1,\dots\,m\}^d}\left(\frac{1}{n} 
   \sum_{i=1}^{n} \varphi_j(X_i)\right) 1_{[t_{j-1}, t_j)}(t).
\end{eqnarray*}
Note that $F_n^{(m)}(t)$ is a piecewise constant approximation to the 
empirical distribution function $F_n(t)$. For $t\in [t_{j-1},t_j)$,
we have the inequality
\[
  F_n(t_{j-2}) \leq F_n^{(m)} (t) \leq F_n(t_{j-1}).
\]
We define further
\[
  F^{(m)}(t) = E \left(F_n^{(m)}(t)\right) = 
 \sum_{j=1}^m E \left(\varphi_j(X_0)\right) 1_{[t_{j-1},t_j)}(t),
\]
and finally the centered and normalized process
\begin{equation}
  U_n^{(m)}(t) = \sqrt{n} \left(F_n^{(m)}(t) - F^{(m)}(t)\right),\quad t\in\R^d.
\label{eq:unm} 
\end{equation}
Theorem~\ref{thm1}  will follow by application of the following Theorem which is proved in Dehling et al. \citeyear{DehDurVol09}.

\begin{thm*}\label{pbil}
Let $ (S,\rho)$ be a complete separable metric space and let 
$X_n$, $X_n^{(m)}$ and $X^{(m)}$, $n,m \geq 1$ be S-valued random variables 
satisfying
\begin{eqnarray}
 && X_n^{(m)}\stackrel {\mathcal{D}}{\longrightarrow } X^{(m)} \mbox{ as }  
 n \rightarrow \infty , \forall m \label{eq:p1}\\ 
 && \lim_{m \rightarrow \infty} \limsup_ {n \rightarrow \infty} 
 P ( \rho (X_n, X_n^{(m)}) \geq \varepsilon) = 0, 
 \forall \varepsilon > 0 .\label{eq:p2}
\end{eqnarray}
Then there exists an S-valued random variable $X$ such that
\[ 
  X_n \stackrel{\mathcal{D}}{\longrightarrow} X \mbox{ as } 
 n \rightarrow \infty .
\]
\end{thm*}

Here we work in the complete separable metric space $(D([-\infty,\infty]^d),d_0)$.
We shall prove separately that (\ref{eq:p1}) and (\ref{eq:p2}) hold for the $D([-\infty,\infty]^d)$-valued random variables $U_n^{(m)}$ in Proposition \ref{pro1} and Proposition \ref{pro2}.

\begin{proposition} \label{pro1}
For any partition $0=r_0 < \ldots < r_m=1$, there exists a 
piecewise constant Gaussian process
$\left(W^{(m)}(t)\right)_{t\in\R^d}$ such that
\[
 \left(U_n^{(m)}(t)\right)_{t\in\R^d} \mathop{\longrightarrow}
 \limits^{\cal{D}} \left(W^{(m)}(t)\right)_{t\in\R^d}.
\]
The sample paths of the processes $\left(W^{(m)}(t)\right)_{t\in\R^d}$ 
are  constant on each of the rectangles
$[t_{j-1}, t_j)$, $1 \leq j \leq m, $ and $W^{(m)}(0) = 0.$
The vector $(W^{(m)} (t_1), \ldots, W^{(m)}(t_m))$
has a multivariate normal distribution with mean zero and covariances
\begin {eqnarray*}
\Cov (W^{(m)}(t_{i-1}), W^{(m)}(t_{j-1}))
&=& \Cov (\varphi_i(X_0), \varphi_j (X_0))\\
&&+ \sum^\infty_{k=1}\Cov (\varphi_i(X_0), \varphi_j (X_k))\\
&&+ \sum^\infty_{k=1}\Cov (\varphi_i(X_k), \varphi_j (X_0))
\end{eqnarray*}
\label{prop:fidi-conv}
\end{proposition}

{\em Proof.}

Use (\ref{eq:clt}) and the Cram\'er-Wold device.
\hfill $\Box$

\begin{proposition} \label{pro2}
For any $\varepsilon,\eta > 0$ there exists a partition $0=r_0<\ldots<r_m=1$ 
such that
\[
  \limsup_{n\rightarrow\infty} 
  P \left(\sup\limits_{t\in\R^d}\left| U_n(t) - U_n^{(m)}(t)\right|>
 \varepsilon\right) \leq \eta.
\]
\label{prop:ep-appr}
\end{proposition}

{\em Proof.}

From here, we assume the partition $0=r_0 < \ldots < r_m=1$ is a regular partition of step $h=m^{-1}$.
Let $j=1,\dots,m$.
On the interval $[r_{j-1}, r_j]$ we introduce a sequence of refining partitions
\[
 r_{j-1} = s_{j,0}^{(k)} < s_{j,1}^{(k)} < \ldots < s^{(k)}_{j,2^k} = r_j
\]
by
\[
 s_{j,l}^{(k)} = r_{j-1} + l \cdot \frac{h}{2^k}\quad,
 \quad 0 \leq l \leq 2^k.
\]
For each $i\in\{1,\dots,d\}$, let us define
$$
s_{i,j,l}^{(k)}=F_i^{-1}(s_{j,l}^{(k)})\quad,\quad 0 \leq l \leq 2^k.
$$
We now have partitions of $[t_{i,j-1},t_{i,j}]$,
\[
 t_{i,j-1}= s_{i,j,0}^{(k)} < s_{i,j,1}^{(k)} < \ldots < s^{(k)}_{i,j,2^k} = t_{i,j}.
\]
For convenience, for $j>1$, we also consider the points
$$
s_{i,j,-1}^{(k)}=F_i^{-1}\left(r_{j-1} - \frac{h}{2^k}\right)=s_{i,j-1,2^k-1}^{(k)}
$$
and for $j<m$, we consider the points
$$
s_{i,j,2^k+1}^{(k)}=F_i^{-1}\left(r_{j-1} + (2^k+1)\frac{h}{2^k}\right)=s_{i,j+1,1}^{(k)}.
$$

For any $t\in [t_{i,j-1}, t_{i,j})$ and $k\geq 0$ we define the index
\[
 l_{i,j}(k,t) = \max \left\{l: s_{i,j,l}^{(k)} \leq t \right\}.
\]

\bigskip

Now fix $j\in\{1,\dots,m\}^d$ (then the index related to $j$ will be forgotten). 
For $l=(l_1,\dots,l_d)$ and $t\in[t_{j-1},t_j)$, we write
\[
 s_l^{(k)}=(s_{1,j_1,l_1}^{(k)},\ldots,s_{d,j_d,l_d}^{(k)})
\]
and
\[
 l(k,t) = (l_{1,j_1}(k,t_1),\ldots,l_{d,j_d}(k,t_d)).
\]

In this way we obtain a chain, 
\[
  t_{j-1} = s_{l(0,t)}^{(0)} \leq s_{l(1,t)}^{(1)} \leq \ldots
  \leq s_{l(k,t)}^{(k)} \leq t \leq s_{l(k,t)+1}^{(k)},
\]
linking the point $t_{j-1}$ to $t$. 
Note that for 
$t\in [t_{j-1}, t_j)$ we have by definition
$U_n^{(m)}(t) = U_n^{(m)}(t_{j-1})$.

We define the functions $\psi^{(k)}_l$,
$k \geq 0$, $ l\in\{0,\dots,2^k+1\}^d$, in the following way :
We first define, for $i=1,\dots,d$ and $l\in\{0,\dots,2^k+1\}$,
\[
\psi^{(k)}_{i,l}(x_i) =
\left\{
\begin{array}{ll}
 0&\mbox{ if } j_i=1\mbox{ and } l=0\\
1 &\mbox{ if } j_i=m\mbox{ and } l\ge 2^k\\
 \displaystyle\varphi\left(\frac{x_i-s_{i,j_i,l}^{(k)}}{|s_{i,j_i,l}^{(k)}-s_{i,j_i,l-1}^{(k)}|}\right)&\mbox{ otherwise}
\end{array}
\right.
\]
where $\varphi$ is defined as in (\ref{eq:phi})  (eventually, we use $\frac{1}{\infty}=0$).
Then we set, for $ l=(l_1,\dots,l_d)\in\{0,\dots,2^k+1\}^d$,
\[
 \psi^{(k)}_l(x) = \prod_{i=1}^{d}\psi^{(k)}_{i,l_i}(x).
\]
Observe that by definition of $s^{(k)}_{l(k,t)}$ and of $\psi^{(k)}$, 
$\psi^{(0)}_{l(0,t)}(x)=\varphi_j(x)$ and 
\[
 \varphi_j(x)\le\psi^{(1)}_{l(1,t)}(x)\le\dots\le\psi^{(k)}_{l(k,t)}(x)\le1_{(-\infty,t]} (x)\le\psi^{(k)}_{l(k,t)+2} (x).
\]

In this way we get
\begin{eqnarray}
 F_n(t) - F_n^{(m)}(t) 
  &=& \sum^K_{k=1} \frac{1}{n} 
   \sum^n_{i=1} \left(\psi^{(k)}_{l(k,t)}(X_i)
       -\psi^{(k-1)}_{l(k-1,t)}(X_i)\right)\nonumber \\
  && + \frac{1}{n}\sum^n_{i=1}\left( 1_{(-\infty, t]}(X_i)
    - \psi^{(K)}_{l(K,t)}(X_i)\right)
\label{eq:fn-fnm}
\end{eqnarray}
where $K$ is some integer to be chosen later.

From (\ref{eq:fn-fnm}) we get by centering and normalization
\begin{eqnarray*}
 U_n(t) - U_n^{(m)}(t)
 &=& \sum^K_{k=1} 
 \frac{1}{\sqrt{n}}\sum^n_{i=1}
  \left\{ 
 \left( \psi^{(k)}_{l(k,t)}( X_i)-E \psi^{(k)}_{l(k,t)}(X_i)
 \right) \right.  \\
 &&-\left. \left(\psi^{(k-1)}_{l(k-1,t)}(X_i)
  -E\psi^{(k-1)}_{l(k-1,t)}(X_i)\right)\right\} \\
 && + \frac{1}{\sqrt{n}}  \sum^{n}_{i=1}
  \left\{\left(1_{(-\infty,t]}(X_i)-F(t)\right)\right.\\
 && - \left. \left(\psi^{(K)}_{l(K,t)}(X_i) 
 -E \psi^{(K)}_{l(K,t)}(X_i)\right)\right\}.
\end{eqnarray*}
For the last term on the r.h.s. we have the following upper and lower bounds,
\begin{eqnarray*}
 &&
 \frac{1}{\sqrt{n}} \sum^n_{i=1} \left\{ 
 \left(1_{(-\infty, t]} (X_i)-F(t)\right) - 
 \left(\psi^{(K)}_{l(K,t)} (X_i) - E \psi^{(K)}_{l(K,t)}
 (X_i)\right)\right\} \\
&& \le  \frac{1}{\sqrt{n}} \sum^n_{i=1} \left\{  
 \left( \psi^{(K)}_{l(K,t)+2} (X_i) - E \psi^{(K)}_{l(K,t)+2} (X_i) \right)- \left( \psi^{(K)}_{l(K,t)} (X_i) 
 - E \psi^{(K)}_{l(K,t)}(X_i) \right) \right\}\\
&&\quad +\sqrt{n} \left( E\psi^{(K)}_{l(K,t)+2} (X_i)-F(t)\right)
\end{eqnarray*}
and
\begin{eqnarray*}
 && \frac{1}{\sqrt{n}}\sum^{n}_{i=1}  \left\{ 
 \left( 1_{(-\infty,t]}(X_i) - F(t)\right)-
 \left(\psi^{(K)}_{l(K,t)}(X_i) - E \varphi^{(K)}_{l(K,t)} 
 (X_i)\right)\right\} \\
 && \geq - \sqrt{n}\left(F(t) - E \psi^{(K)}_{l(K,t)} 
  (X_i)\right).
\end{eqnarray*}
Now choose $K = 4 + \left\lfloor \log\left(d \frac{\sqrt{n}h}{\varepsilon}\right)\log^{-1}(2)\right\rfloor$
 and note that 
\begin{equation*}
\frac{\varepsilon}{2^4}\leq  d\sqrt {n} \frac{h}{2^K} \leq \frac{\varepsilon}{2^3}.
\end{equation*}
We thus have
\begin{eqnarray*}
\sqrt{n}\left| E \psi^{(K)}_{l(K,t)+2}(X_i)- 
E \psi^{(K)}_{l(K,t)}(X_i)\right|
&\leq&\sqrt{n}\sum_{i=1}^d\left|F_i(s^{(K)}_{l_{i,j_i}(K,t)+2})-F_i(s^{(K)}_{l_{i,j_i}(K,t)-1})\right|\\
&\leq&\frac{\varepsilon}{2}.
\end{eqnarray*}
Therefore, since
\[
 \psi^{(K)}_{l(K,t)} ( X_i) \leq  1_{(-\infty,t]} (X_i) \leq  \psi^{(K)}_{l(K,t)+2} (X_i),
\]
we get for all $t \in [t_{j-1}, t_j]$, 
\begin{eqnarray*}
\left| U_n(t)-U_n^{(m)} (t)\right|
 &\leq& \sum^K_{k=1} \frac{1}{\sqrt{n}} 
  \left| \sum^n_{i=1}
  \left\{  \left( \psi^{(k)}_{l(k,t)} (X_i) - 
  E \psi^{(k)}_{l(k,t)} (X_i) \right) \right.\right. \\
&& \left.\left.\quad - \left( \psi^{(k-1)}_{l(k-1,t)}  (X_i) - 
  E \psi^{(k-1)}_{l(k-1,t)}  (X_i)\right)  \right\} \right| 
  \\
&& + \frac{1}{\sqrt{n}} \left|  \sum\limits^n_{i=1} \left\{  
 \left( \psi^{(K)}_{l(K,t)+2}  (X_i) - E \psi^{(K)}_{l(K,t)+2}  (X_i)\right)\right.\right. \\
&&\left.\left. \quad -\left( \psi^{(K)}_{l(K,t)}  (X_i) - E \psi^{(K)}_{l(K,t)} 
  (X_i) \right)   \right\} 
   \right| \\
&& \quad + \frac{\varepsilon}{2}.
\end{eqnarray*}
Note that by definition of $l(k,t)$ and of $s_{l}^{(k)}$, we have
\[
l(k-1,t)=\left\lfloor \frac{l(k,t)}{2} \right\rfloor
\]
where the integer part $\lfloor .\rfloor$ is taken on each coordinate.
We infer
\begin{eqnarray*}
\sup_{t_{j-1} \le t \le t_j}  \left| 
    U_n(t) - U^{(m)}_n (t) \right|
&\le&  \sum^K_{k=1} \frac{1}{\sqrt{n}} \max_{l\in\{0,\dots,2^k\}^d}
   \left|  \sum^n_{i=1} 
  \left( ( \psi^{(k)}_l (X_i) - E \psi^{(k)}_l
  (X_i) )\right. \right. \\
&& \qquad \quad \left. \left. - ( \psi^{(k-1)}_{\lfloor\frac{l}{2}\rfloor} (X_i) - E \psi^{(k-1)}_{\lfloor\frac{l}{2}\rfloor}
  (X_i)) \right)  \right|\\
&& \quad + \frac{1}{\sqrt{n}} \max_{l\in\{0,\dots,2^K\}^d}   \left|  
  \sum^n_{i=1} \left( (\psi^{(K)}_{l+2} (X_i) - E \psi^{(K)}_{l+2} 
  (X_i)) \right. \right. \\
&& \left.\left. \qquad \quad - ( \psi^{(K)}_l (X_i) - E \psi^{(K)}_l
  (X_i)) \right) \right| \\
&& \quad + \frac{\varepsilon}{2}.
\end{eqnarray*}

\medskip

Now, taking $\varepsilon_k   = \frac{\varepsilon}{4 k (k+1)} $, we obtain
\begin{eqnarray*}
&&\hspace{-15mm} P \left(\sup_{t_{j-1}\le t \le t_j} 
  \left|U_n (t) - U_n^{(m)} (t) \right| \ge\varepsilon\right)\\
& \le &\sum\limits^K_{k=1} \sum\limits_{l\in\{0,\dots,2^k\}^d} 
P \left( \frac{1}{\sqrt{n}}\right. \left| \sum\limits^n_{i=1} \right. \left\{
\left( \psi^{(k)}_l (X_i) - E \psi^{(k)}_l 
(X_i) \right) \right. \\
&& \quad \left.\left.\left.-\left( \psi^{(k-1)}_{\lfloor\frac{l}{2}\rfloor} (X_i) - E \psi^{(k-1)}_{\lfloor\frac{l}{2}\rfloor} 
 (X_i)
\right) \right\} \right| \ge \varepsilon_k \right) \\
&& + \sum_{l\in\{0,\dots,2^K\}^d} P \left( \frac{1}{\sqrt{n}}\right.  \left| 
\sum^n_{i=1}\left\{ \left( \psi^{(K)}_{l+2} (X_i) - 
E \psi^{(K)}_{l+2} (X_i) \right)\right.\right. \\
&& \quad \left.\left.\left.-\left(\psi^{(K)}_l(X_i) - E \psi^{(K)}_l (X_i)
\right)  \right\} \right|  \ge \frac{\varepsilon}{4}\right).
\end{eqnarray*}
At this point we shall use Markov's inequality at the order $2p$ together with the $2p$-th moment 
bound (\ref{eq:moment}) for an integer $p$ such that
\begin{equation}\label{gamma}
p>d\frac{r\gamma}{\gamma-r}.
\end{equation}

First, remark that these following bounds also hold. Since $\psi^{(k)}_{i,l_i} - \psi^{(k-1)}_{i,\lfloor\frac{l_i}{2}\rfloor}$ vanishes outside $[s_{i,j_i,{\lfloor\frac{l_i}{2}\rfloor}-1}^{(k-1)},s_{i,j_i,l_i}^{(k)}]$, we have
\begin{eqnarray*}
 \left\| \psi^{(k)}_l(X_0) - \psi^{(k-1)}_{\lfloor\frac{l}{2}\rfloor} 
 (X_0) \right\|_r 
&\leq & \left(\sum_{i=1}^d\left| F_i(s_{i,j_i,l_i}^{(k)}) - F_i(s_{i,j_i,{\lfloor\frac{l_i}{2}\rfloor}-1}^{(k-1)}) \right|\right)^{\frac{1}{r}}\\
&\leq & d \max_{i=1}^{d}\left| F_i(s_{i,j_i,l_i}^{(k)}) - F_i(s_{i,j_i,l_i-3}^{(k)}) \right|^{\frac{1}{r}}\\
& =&  \left(\frac{3dh}{2^k}\right)^{\frac{1}{r}}
\end{eqnarray*}
and in the same way
\begin{eqnarray*}
 \left\| \psi^{(K)}_{l+2}(X_0) - \psi^{(K)}_{l} 
 (X_0) \right\|_r&\le & \left(\frac{3dh}{2^K}\right)^{\frac{1}{r}}.
\end{eqnarray*}

\medskip

Now by (\ref{eq-modulus}), if $k$ is big enough, we have
\begin{eqnarray*}
 \left\| \psi^{(k)}_l \right\|
&\le&1+d\max_{i=1}^d\frac{1}{|s_{i,j_i,l_i}^{(k)}-s_{i,j_i,l_i-1}^{(k)}|^\alpha}\\
&\leq& 1+d\left[\inf\left\{s>0 : \exists i\in\{1,\dots,d\},\forall t, F_i(t+s)-F_i(t)\geq \frac{h}{2^k}\right\}\right]^{-\alpha}\\
&\leq& 1+d\left[\inf\left\{s>0 : D|\log(s)|^{-\gamma}\geq \frac{h}{2^k}\right\}\right]^{-\alpha}\\
&=& 1+ d \exp\left(\alpha\left(\frac{D2^k}{h}\right)^\frac{1}{\gamma}\right).
\end{eqnarray*}
Thus, there is a positive constant $B$ such that for arbitrary $k\ge 0$, 
$$
\left\| \psi^{(k)}_l \right\|\le B\exp\left(\alpha\left(\frac{D2^k}{h}\right)^\frac{1}{\gamma}\right).
$$

Therefore, applying successively Markov's inequality at the order $2p$, the $2p$-th moment 
bound (\ref{eq:moment}) and the preceding inequalities, we get
\begin{eqnarray*}
  && \hspace*{-20mm} P \left(\sup\limits_{ t_{j-1} \leq t \leq t_j}
 \left| U_n (t) - U_n (t_j)\right|\geq \varepsilon\right) \\
 &\leq&  C \sum_{i=1}^{p}
\sum\limits^K_{k=1} 2^{dk} \frac{(k(k+1))^{2p}}{\varepsilon^{2p}} 
 \frac{1}{n^p}n^i\left(\frac{h^i}{2^{ik}}\right)^{\frac{1}{r}}\log^{2p-i} \left(a+B\exp\left(\alpha\left(\frac{D2^k}{h}\right)^\frac{1}{\gamma}\right)
  \right)\\
&\leq&  C \sum_{i=1}^{p-1}\frac{h^d}{\varepsilon^{2p}n^{p-i}} 
\sum\limits^K_{k=1} \left(\frac{2^k}{h}\right)^{d-\frac{i}{r}} 
 k^{4p}\left(\frac{2^k}{h}\right)^\frac{2p-i}{\gamma}\\
&&+ C\frac{h^d}{\varepsilon^{2p}}\sum\limits^K_{k=1}
\left(\frac{2^k}{h}\right)^{d-\frac{p}{r}}k^{4p}\left(\frac{2^k}{h}\right)^\frac{p}{\gamma}\\
&\leq&  C \sum_{i=1}^{p-1}\frac{h^d}{\varepsilon^{2p+d-i+\frac{2p-i}{\gamma}}} 
 (\sqrt{n})^{d-\frac{i}{r}+\frac{2p-i}{\gamma}-2(p-i)} 
 K^{4p+1}\\
&&+C\frac{h^{\frac{p}{r}-\frac{p}{\gamma}}}{\varepsilon^{2p}}\sum\limits^\infty_{k=1}
2^{(d-\frac{p}{r}+\frac{p}{\gamma})k}k^{4p},
\end{eqnarray*}
where $C$ always denotes a positive constant, but its value changes from line to line.

By Condition (\ref{gamma}), the series $\sum\limits^\infty_{k=1}
2^{(d-\frac{p}{r}+\frac{p}{\gamma})k}k^{4p}$ converges and there exists an $A>0$ such that
for all $i\in\{1,\dots,p-1\}$, 
$$
d-\frac{i}{r}+\frac{2p-i}{\gamma}-2(p-i)<-A.
$$
Finally, using $mh=1$, we have
\begin{eqnarray*}
&& \hspace*{-20mm} P \left( \sup_{0\leq t \leq 1} 
 \left| U_n (t) - U_n^{(m)}(t) \right| 
\geq \varepsilon \right) \\
&\leq& \sum\limits_{j\in\{1,\dots,m\}^d} P \left( \sup_{t_{j-1} \leq t \leq t_j} 
 \left| U_n(t) - U_n^{(m)} (t) \right| \geq \varepsilon  \right)\\
 &\leq&  C (p-2)\frac{1}{\varepsilon^{4p+d}} 
 (\sqrt{n})^{-A}\left(4+\log\left(d\frac{\sqrt{n}h}{\varepsilon}\right)\right)^{4p+1}
+C\frac{h^{\frac{p}{r}-d-\frac{p}{\gamma}}}{\varepsilon^{2p}}.
\end{eqnarray*}
The first summand converges to zero as 
$n \rightarrow \infty$ and, since $\frac{p}{r}-d-\frac{p}{\gamma}>0$, the second can be made arbitrarily  
small by choosing a partition that is fine enough (i.e. $h$ small). \hfill $\Box$

\begin{rem}
 The tightness of the empirical process can be proved using exactly the same proof than the one of Proposition 2.3 in 
Dehling et al. \citeyear{DehDurVol09}. The almost sure continuity of the limit process follows.
\end{rem}

\medskip

\begin{rem}
Assume that conditions (\ref{eq:clt}) and (\ref{eq:moment}) hold for a space of $C^k$-functions instead of a space of H\"older continuous functions. It is clear that the same technique works (in the proof, take for example $\varphi(x)=1_{(-\infty,-1]}(x)-\sin\left(\pi x +\frac{\pi}{2}\right)1_{(-1,0]}(x)$).
\end{rem}

\medskip

\section{Moment bounds for partial sums}
\label{sec:mb}

For a function $\varphi:\R^d\longrightarrow\R$, we consider the partial sum 
$$
S_n\left(\varphi\right)=\sum_{i=0}^{n-1}\varphi\left(X_i\right).
$$

Multiple mixing properties allow us to obtain some useful moment inequalities. In this section, we show how the multiple mixing property implies the 2p-th moment bound which is required in Theorem \ref{thm1}.

\begin{thm}\label{MB}
Let $\left(X_n\right)_{n\ge 0}$ be a stationary process having a multiple mixing property on $\B$ and 
$\varphi\in\B$ such that $\esp_\nu\left(\varphi\right)=0$ and $\sup_x|\varphi(x)|\le 1$. Then for all $p\ge 1$,
$$
\esp_\nu\left[S_n\left(\varphi\right)^{2p}\right]\leq
K\sum_{i=1}^p n^i\|\varphi\left(X_0\right)\|_r^i\log^{2p-i}(\|\varphi\|+\theta^{-1})
$$
and
$$
\left|\esp_\nu\left[S_n\left(\varphi\right)^{2p+1}\right]\right|\leq
K\sum_{i=1}^p n^i\|\varphi\left(X_0\right)\|_r^i\log^{2p-i+1}(\|\varphi\|+\theta^{-1})
$$
where $K$ is a constant which does not depend on $n$ or $\varphi$ and $r\ge 1$ is given by (\ref{mm}).
\end{thm}

\medskip

\paragraph*{Proof of Theorem~\ref{MB}.}

Let us consider the assumptions of Theorem~\ref{MB} hold and 
let $\varphi\in\B$ with $\esp\left(\varphi\right)=0$ and $\sup_x|\varphi(x)|\le 1$ be fixed.
We use the notation $a^*_n=\sum_{i=1}^n a_i$.
\medskip

\noindent\textbf{Notations:}
For all $p\ge 1$, we define
$$
I_n(p)=\sum_{\begin{array}{c}0\le i_1,\dots,i_p\le n-1\\i_p^*\le n-1\end{array}}
|\esp(\varphi(X_0)\varphi(X_{i_1^*})\dots\varphi(X_{i_p^*}))|
$$
and $I_n(0)=|\esp(\varphi(X_0))|=0$.

\medskip

As the process is stationary, for $p\ge 1$, we have
\begin{equation*}
\left|\esp\left[S_n\left(\varphi\right)^p\right]\right|\leq p!nI_n(p-1).
\end{equation*}
So, to prove Theorem~\ref{MB}, it is sufficient to prove the following lemma.

\begin{lem}\label{key3}
For all $p\ge 1$,
$$
I_n(2p-1)\leq
K\sum_{i=1}^p n^{i-1}\|\varphi\left(X_0\right)\|_r^i\log^{2p-i}(\|\varphi\|+\theta^{-1})
$$
and
$$
I_n(2p)\leq
K\sum_{i=1}^{p} n^{i-1}\|\varphi\left(X_0\right)\|_r^i\log^{2p-i+1}(\|\varphi\|+\theta^{-1})
$$
where $K$ is a constant which does not depend on $n$ or $\varphi$.
\end{lem}

\medskip

\noindent\textbf{Notations:}
For all $p\ge 1$ and $q\in\{1,\dots,p\}$, we define
$$
J_n(p,q)=\sum_{i_q=0}^{n-1}\sum_{\begin{array}{c}0\le i_1,\dots,i_{q-1},i_{q+1},\dots,i_p\le i_q\\i_p^*\le n-1\end{array}}|\esp(\varphi(X_0)\varphi(X_{i_1^*})\dots\varphi(X_{i_p^*}))|.
$$

We have
$$
I_n(p)\le \sum_{q=1}^p J_n(p,q).
$$

To prove Lemma \ref{key3}, we will use the following lemma.

\begin{lem}\label{key2}
 For all $p\in\N^*$ and $q\in\{1,\dots,p\}$,
$$
J_n(p,q)\le C\|\varphi(X_0)\|_r\log^p(\|\varphi\|+\theta^{-1})+nI_n(q-1)I_n(p-q),
$$
where $C$ is a constant which does not depend on $n$ or $\varphi$.
\end{lem}

\medskip

{\em Proof of Lemma \ref{key2}.}

Let $n_0$ be a positive integer such that 
$$
\frac{\log(\|\varphi\|+\theta^{-1})}{-\log\theta}< n_0\leq \frac{\log(\|\varphi\|+\theta^{-1})}{-\log\theta}+1.
$$
We thus have the inequality  $\theta^{n_0}\|\varphi\|\leq 1$ and $n_0\ge 2$. 

\medskip

We have
\begin{equation*}
J_n(p,q)\le \sum_{i_q=0}^{n-1}\sum_{\begin{array}{c}0\le i_1,\dots,i_{q-1},i_{q+1},\dots,i_p\le i_q\\i_p^*\le n-1\end{array}}
\left[A_{i_1,\dots,i_p}+B_{i_1,\dots,i_p}\right]
\end{equation*}
where
$$A_{i_1,\dots,i_p}=|\Cov(\varphi(X_0)\varphi(X_{i_1^*})\dots\varphi(X_{i_{q-1}^*}),
\varphi(X_{i_q^*})\varphi(X_{i_{q+1}^*})\dots\varphi(X_{i_p^*})|
$$
and
$$ B_{i_1,\dots,i_p}=|\esp(\varphi(X_0)\varphi(X_{i_1^*})\dots\varphi(X_{i_{q-1}^*}))||\esp(\varphi(X_0)\varphi(X_{i_{q+1}})\dots\varphi(X_{i_p^*-i_q^*}))|.
$$

\medskip

Using H\"older inequality for $i_q=0$ to $n_0\ell-2$ and multiple mixing property (\ref{mm}) for $i_q\ge n_0\ell$
(where $\ell$ comes from (\ref{mm})), we obtain
\begin{eqnarray*}
&&\hspace{-10pt}
\sum_{i_q=0}^{n-1}\sum_{\begin{array}{c}0\le i_1,\dots,i_{q-1},i_{q+1},\dots,i_p\le i_q\\i_p^*\le n-1\end{array}}A_{i_1,\dots,i_p}\\
&&\le C\sum_{i_q=0}^{n_0\ell-2}(i_q+1)^{p-1}\|\varphi(X_0)^{p+1}\|_1 +C\sum_{i_q=n_0\ell-1}^{n}(i_q+1)^{p-1}\kappa\|\varphi(X_0)\|_r\theta^{i_q}\|\varphi\|^{\ell}Q(i_q)\\
&&\le C(n_0-1)^p\|\varphi(X_0)\|_r
+C\kappa\|\varphi(X_0)\|_r\sum_{i_q=n_0\ell-1}^{\infty}(i_q+1)^{p-1}\theta^{i_q-n_0\ell}Q(i_q)
\end{eqnarray*}
where $Q(i_q)=\sum_{0\le i_1,\dots,i_{q-1},i_{q+1},\dots,i_p\le i_q}P(i_1,\dots,i_p)$ is a polynomial function of $i_q$.
Then $\sum_{i_q=0}^{\infty}\theta^{i_q}Q(i_q)$ converges and we deduce that
$$
\sum_{i_q=n_0\ell-1}^{\infty}(i_q+1)^{p-1}\theta^{i_q-n_0\ell}Q_n(i_q)\leq 
C(n_0-1)^{p-1}
$$
where $C$ is independent of $n_0$.
Thus, since $n_0-1\le C\log(\|\varphi\|+\theta^{-1})$,
$$
\sum_{i_q=0}^{n-1}\sum_{\begin{array}{c}0\le i_1,\dots,i_{q-1},i_{q+1},\dots,i_p\le i_q\\i_p^*\le n-1\end{array}}A_{i_1,\dots,i_p}
\le C\kappa\|\varphi(X_0)\|_r\log^p(\|\varphi\|+\theta^{-1}).
$$
On the other hand,
$$
\sum_{\begin{array}{c}0\le i_1,\dots,i_{q-1},i_{q+1},\dots,i_p\le i_q\\i_p^*\le n-1\end{array}}B_{i_1,\dots,i_p}
\le I_n(q-1)I_n(p-q).
$$
Therefore,
\begin{equation*}
 J_n(p,q)\le C\kappa\|\varphi(X_0)\|_r\log^p(\|\varphi\|+\theta^{-1})+nI_n(q-1)I_n(p-q).
\end{equation*}
\cqfd

\medskip

{\em Proof of Lemma \ref{key3}.}

We proceed by induction.
We have, $I_n(1)=J_n(1,1)$. Then, by Lemma \ref{key2},
$$
I_n(1)\le C\|\varphi(X_0)\|_r\log(\|\varphi\|+\theta^{-1}).
$$

In the same way, $I_n(2)\le J_n(2,1) + J_n(2,2)$. Then, by Lemma \ref{key2},
$$
I_n(2)\le C\|\varphi(X_0)\|_r\log^2(\|\varphi\|+\theta^{-1}).
$$

In the general case, by Lemma \ref{key2},
\begin{eqnarray*}
I_n(p)&\le& \sum_{q=1}^p J_n(p,q)\\
&\le&Cp\|\varphi(X_0)\|_r\log^p(\|\varphi\|+\theta^{-1})+n\sum_{q=2}^{p-1}I_n(q-1)I_n(p-q).
\end{eqnarray*}
Studying $\sum_{q=2}^{p-1} I_n(q-1)I_n(p-q)$ according to the parity of $p$, we deduce the inequalities of Lemma \ref{key3}.
\cqfd

\medskip

\section{Markov chains and dynamical systems with a spectral gap}\label{spectral}

Let $\left(X_n\right)_{n\ge 0}$ be a homogeneous Markov chain with a stationary measure $\nu$.
Denote by $P$ the associated Markov operator and $E$ the state space.
Consider a Banach algebra $\left(\B,\|. \|\right)$ of $\nu$-measurable functions from $E$ to $\R$, which contains
the function $\bf{1}=\ind_E$ and which is continuously included in $\left(\L^s\left(\nu\right), \|.\|_s\right)$
for some $s\in[1,+\infty]$,\\
i.e. $\exists C>0$ such that $\forall f\in\B,$ 
\begin{equation}\label{cont}
\|f\|_s\leq C\|f\|.
\end{equation}

We say that the Markov chain $\left(X_n\right)_{n\ge 0}$ is {\it $\B$-geometrically ergodic} or {\it strongly ergodic} (with respect to $\B$)
if there exist $\kappa>0$ and $0<\theta<1$ such that for all $f\in\B$, 
\begin{equation}\label{exp}
\|P^nf-\Pi f\|\leq \kappa\theta^n\|f\|
\end{equation}
where $\Pi f=\esp_\nu\left(f\right)\bf{1}$.

\medskip

Strong ergodicity  corresponds to the fact that the Markov
transition operator acting on $\B$ has $1$ as simple eigenvalue and the rest of the spectrum is included in a closed ball of radius strictly smaller than $1$ \citeaffixed{HenHer01}{see}.

\begin{lem}\label{key1}
Let $\left(X_n\right)_{n\ge 0}$ be a $\B$-geometrically ergodic Markov chain, then it satisfies the multiple mixing property (\ref{mm})on $\B$ with $r=\frac{s}{s-1}$ ($r=1$ if $s=\infty$).
\end{lem}

\medskip

{\em Proof of Lemma \ref{key1}.}

Let $\F_i$ be the $\sigma$-algebra generated by the $X_j$, $j\le i$ and let $\varphi$ belongs to $\B$ such that $\esp_\nu\left(\varphi\right)=0$ and $\|\varphi\|_\infty\le 1$.
Using the operator properties, we have
\begin{eqnarray*}
&&\hspace{-10pt}|\Cov(\varphi(X_0)\varphi(X_{i_1^*})\dots\varphi(X_{i_{q-1}^*}),
\varphi(X_{i_q^*})\varphi(X_{i_{q+1}^*})\dots\varphi(X_{i_p^*})|\\
&&= |\esp\left[\varphi(X_0)\varphi(X_{i_1^*})\dots\varphi(X_{i_{q-1}^*})\left[
\esp(\varphi(X_{i_q^*})\esp(\varphi(X_{i_{q+1}^*})\dots\esp(\varphi(X_{i_p^*})|\F_{i_{p-1}^*})\dots|\F_{i_{q}^*})|\F_{i_{q-1}^*})
\right.\right.\\
&&\hspace{250pt}\left.\left.-\esp(\varphi(X_{i_q^*})\varphi(X_{i_{q+1}^*})\dots\varphi(X_{i_p^*}))\right]\right]|\\
&&\le \|\varphi(X_0)\varphi(X_{i_1^*})\dots\varphi(X_{i_{q-1}^*})\|_r\\
&&\hspace{100pt}\qquad\|P^{i_q}(\varphi P^{i_{q+1}}(\varphi\dots P^{i_p}\varphi))(X_0)-\pi(\varphi P^{i_{q+1}}(\varphi\dots P^{i_p}\varphi))(X_0)\|_s\\
&&\le C\kappa\|\varphi(X_0)^q\|_r\theta^{i_q}\|\varphi P^{i_{q+1}}(\varphi\dots P^{i_p}\varphi)\|.
\end{eqnarray*}
Further,
$$
\|\varphi P^{i_{q+1}}(\varphi\dots P^{i_p}\varphi)\|\le C\|\varphi\|^{p-q+1}
$$
and the result follows.
\cqfd

The corresponding result holds in the setting of dynamical systems.
Let $(\Omega,\A,\mu)$ be a probability space and $T$ a measurable measure preserving transformation.
Let us consider the Perron-Frobenius operator (or the transfer operator) of $T$, $P:\L^1(\mu)\longrightarrow\L^1(\mu)$
defined by
$$
\int_{\Omega}Pf(x)g(x)\dif\mu(x)=\int_{\Omega}f(x)g\circ T(x)\dif\mu(x)
$$
for all $f\in\L^1(\mu)$ and $g\in\L^\infty(\mu)$.

We assume there exists a Banach algebra $(\B,\|.\|)$ of $\mu$-measurable functions from $\Omega$ to $\R$ which
contains $\bf{1}$ and satisfies (\ref{cont}) and that $P$ verifies:\\
there exist $\kappa>0$ and $0<\theta<1$ such that for all $f\in\B$, 
$$
\|P^nf-\Pi f\|\leq \kappa\theta^n\|f\|
$$
where $\Pi f=\esp_\mu\left(f\right)\bf{1}$.

\begin{lem}
Let $f\in\B$ and $X_i=f\circ T^i$, $i\le 0$. Then $(X_i)_{i\le 0}$ satisfies the multiple mixing property (\ref{mm}) on $\B$ with $r=\frac{s}{s-1}$ ($r=1$ if $s=\infty$).
\end{lem}

\medskip

For both setting, if the space $\B$ is enough well adapted to approximate indicator functions, then Theorem~\ref{thm1} applies.
Examples of Section 4 in Dehling et al. \citeyear{DehDurVol09} can be generalized in higher dimensions.
In particular, we can state a result for multidimensional linear processes.

\medskip

\paragraph*{Linear processes.}

Let $(\xi_i)_{i\ge Z}$ be a sequence of i.i.d. random variables in $\R^d$ with $\|\xi_0\|_\infty<\infty$.
Let $(a_i)_{i\ge 0}$ be a sequence of endomorphisms of $\R^d$ such that $\|a_i\|_\infty\le \theta^i$ (for a $\theta<1$).
Define the linear process $X_k=\sum_{i\ge 0}a_i(\xi_{k-i})$ and assume that its distribution function $F$ satisfies
condition (\ref{eq-modulus}). Then we can show that the process $(X_k)_{k\ge 0}$ satisfies condition (\ref{exp}) for the space of bounded Lipschitz functions on $\R^d$. We deduce conditions (\ref{eq:clt}) and (\ref{eq:moment}) and Theorem~\ref{thm1} leads to the following corollary.
\begin{cor}
 If the distribution function of $X_0$ verifies condition (\ref{eq-modulus}) with $r=1$, then the multivariate empirical process associated to $(X_k)_{k\ge 0}$ converges in distribution to an almost surely continuous Gaussian process.
\end{cor}
This result was already proved by Dedecker \citeyear{Ded09}, using a different technique.

\medskip

\paragraph*{Random iterative Lipschitz models.}

Here, let us focus on application concerning random iterative Lipschitz models.
This example has been investigated before by Wu and Shao \citeyear{WuSha04} or Dedecker and Prieur \citeyear{DedPri07}, with different techniques and under different conditions. Here we want to show that our technique also applies in this situation.

\medskip

Let $g:\R^d\times \R \longrightarrow \R^d$ be a measurable function and let $(Y_n)_{n\ge 1}$ be an $\R$-valued i.i.d. process. Let $X_0$ be an $\R^d$-valued random variable independent of $(Y_n)_{n\ge 1}$.
Define the Markov chain $(X_n)_{n\in\N}$ by
$$
X_n = g(X_{n-1},Y_n),\quad n\ge 1.
$$

\medskip

Assume that for all $y\in\R$, $g(.,y)$ is Lipschitz. Define the Lipschitz constant
\[
 K(y):=\sup_{x,x'\in E, x\neq x'}\frac{|g(x,y)-g(x',y)|}{|x-x'|}
\]
and suppose that there exists $\gamma_0>1$ such that 
\[
 E[(1+K(Y_1)+|g(0,Y_1)|)^{\gamma_0+1}(1+K(Y))]<\infty
\]
and
\[
 E[K(Y_1)\max\{K(Y_1),1\}^{2\gamma_0}]<1.
\]
Let $\gamma\in(0,\gamma_0]$ and consider the Banach space $\B_\gamma$ of functions from $\R^d$ to $\R$
satisfying
\[
 m_\gamma(f)=\sup_{x\neq y}\frac{|f(x)-f(y)|}{|x-y|p(x)^\gamma p(y)^\gamma}<\infty
\]
where $p(x)=1+|x|$. The associated norm is $\|.\|_\gamma=N_\gamma(.)+m_\gamma(.)$ with $N_\gamma(f)=\sup_x\frac{|f(x)|}{p(x)^{\gamma+1}}$.
According to Hennion and Herv\'e \citeyear{HenHer04} Theorem~5.5, the Markov chain has an invariant probability measure $\nu$ such that $\nu(|.|^{\gamma_0+1})<\infty$ and for all $\gamma\in(0,\gamma_0]$, the chain is $\B_\gamma$ geometrically ergodic. 
Further, if $f\in\B_\gamma$, $\nu(|f|^\gamma_0)^\frac{1}{\gamma_0}\le N_\gamma(f)\nu(p^{\gamma_0(\gamma+1)})^\frac{1}{\gamma_0}$. Then for $\gamma_1=\frac{1}{\gamma_0}$, $\B_{\gamma_1}$ is continuously included in $L^{\gamma_0}(\nu)$.
By Lemma \ref{key1}, $(X_n)_{n\ge 0}$ satisfies a multiple mixing property on $\B_{\gamma_1}$, for $r=\frac{\gamma_0}{\gamma_0-1}$.
Note that if $f$ is a bounded Lipschitz function on $\R^d$ then $f\in\B_{\gamma_1}$ and 
$\|f\|_{\gamma_1}\le\|f\|$ (where $\|f\|=\sup|f|+m_0(f)$). Thus we have the multiple mixing property for the space of bounded Lipschitz functions. As a corollary of Theorem~\ref{thm1} we get:

\begin{cor}
Assume that the distribution function of $\nu$ satisfies (\ref{eq-modulus}) for $r=\frac{\gamma_0}{\gamma_0-1}$. If the distribution of $X_0$ is $\nu$, then the empirical process associated to $(X_n)_{n\ge 0}$ converges in distribution to an almost surely continuous Gaussian process.
\end{cor}

\medskip

\section{Ergodic torus automorphisms}\label{sec:torus}

Let $T$ be an ergodic automorphism of the torus of dimension $d$ and $\mu$
the Lebesgue measure on $\T^d$. Then we have the following multiple mixing property.
\begin{proposition}
There exist $C>0$, $0<\gamma <1$, for all $m,p\in \N^*$,
for all bounded $\alpha$-H\"older function $\varphi$ ($\alpha\in(0,1]$) with $\| \varphi\|_\infty\le 1$,
for all $k_1\leq...\leq k_m\leq0\leq l_1\leq...\leq l_p$,
for all $n \in\N$,
\begin{equation}
\left\vert \Cov\left(\prod_{j=1}^m\varphi\circ T^{k_j},\prod_{j=1}^p\varphi\circ T^{l_j+n}\right)\right\vert
\leq C\|\varphi\|_1\|\varphi\|_{\H^\alpha}P(k_1,\dots,k_m)\gamma^ n
\end{equation}
where $P(k_1,\dots,k_m)=\sum_{i=1}^m|k_i|^r$ with $r$ the size of the biggest Jordan's block of $T$ restricted to its neutral subspace.
\end{proposition} 

Almost the same proposition appears in Le Borgne and P\`ene \citeyear{LebPen05}. The slightly modification here is that we keep the $L^1$-norm appearing in the upper bound.

{\em Proof. }

Denote by $E^s$, $E^u$ and $E^c$ the $T$-stable subspaces of $\R^d$ corresponding respectively to the stable , the unstable and the central directions of $T$ (where $T$ is identified to its representative matrix). We have
$$
\R^d=E^s\oplus E^u\oplus E^c
$$
and there exists $\lambda>1$ such that for all $n\ge 1$,
$$
|T^nv|\le \lambda^{-n}|v|\mbox{ for all } v\in E^s,
$$
$$
|T^{n}v|\ge \lambda^n |v|\mbox{ for all } v\in E^u,
$$
$$
|T^nv|\le n^r |v| \mbox{ for all } v\in E^c,
$$
where $|.|$ denotes the maximum norm on $\R^d$ and $r$ is the size of the greatest Jordan block of $T$ restricted to the space
$E^c$.
Further $\mu$ can be written as the product measure of $\mu_s$, $\mu_u$ and $\mu_c$.
Set 
$$
B_i(0,\rho)=\{x\in E^i\,/\,|x|\le \rho\}, \; \rho>0,i=s,u,c.
$$

\medskip

Denote by $\|.\|_{\H^\alpha_s}$ (resp. $\|.\|_{\H^\alpha_{u,c}}$) the $\H^\alpha$-norm is the stable direction (resp. unstable-central direction). Following ideas of the proof of Proposition III.3 in Le Borgne \citeyear{Leb99} \citeaffixed{LebPen05}{see also}, one can prove a result concerning the good distribution of the stable leaves in the torus.

\begin{lem}\label{gdsl} (good distribution of stable leaves)
 There exist $\theta<1$ such that for all $\phi\in\H^\alpha$, $x\in\T^d$ and $\rho>0$,
$$
\frac{1}{\mu_s(T^{-n}B_s(0,\rho))}\left|\int_{T^{-n}B_s(0,\rho)}\phi(T^{-n}x+s)ds\right|\le C\|\phi\|_{\H^\alpha_{u,c}}\theta^ n.
$$ 
\end{lem}

\medskip

Let $\A_0$ be a sub-$\sigma$-algebra of the Borelian one for which the atoms are pieces of stable leaves and set
$\A_n=T^{-n}\A_0$. Let $\phi$ and $\psi$ be two $\Cl$-function with zero mean.
\begin{eqnarray*}
 \Cov(\phi,\psi\circ T^n)&=&\Cov(\phi-\esp(\phi|\A_{-\lfloor\frac{n}{2}\rfloor}),\psi\circ T^n)
+\Cov(\esp(\phi|\A_{-\lfloor\frac{n}{2}\rfloor}),\psi\circ T^n)\\
&\le&\|\psi\|_1\|\phi-\esp(\phi|\A_{-\lfloor\frac{n}{2}\rfloor})\|_\infty+\|\phi\|_1\|\esp(\psi|\A_{\lceil\frac{n}{2}\rceil})\|_\infty.
\end{eqnarray*}
But, since the diameter of the atoms of $\A_{-n}$ decreases exponentially fast,
$$
\|\phi-\esp(\phi|\A_k)\|_\infty\le C\|\phi\|_{\H^\alpha_s}\lambda^{-\frac{n}{2}}
$$
and, by Lemma \ref{gdsl},
$$
\|\esp(\psi|\A_{\lceil\frac{n}{2}\rceil})\|_\infty\le C\|\psi\|_{\H^\alpha_{u,c}}\theta^\frac{n}{2}.
$$
Thus, for $\gamma=\max\{\lambda,\theta\}^\frac{1}{2}<1$, we get
\begin{equation}\label{inn}
|\Cov(\phi,\psi\circ T^n)|\le C(\|\psi\|_1\|\phi\|_{\H^\alpha_s}+\|\phi\|_1\|\psi\|_{\H^\alpha_{u,c}})\gamma^n.
\end{equation}

Further, for all $n\ge 0$ and for all $\phi\in\H^\alpha$, we have 
\begin{equation}\label{inn1}
\|\phi\circ T^n\|_{\H^\alpha_s}\le C\|\phi\|_{\H^\alpha_s}.
\end{equation}
Indeed, if $x\in\T^d$ and $v\in E^s$, by linearity of the map $T$,
\begin{eqnarray*}
|\phi\circ T^n(x)-\phi\circ T^n(x+v)|&=&|\phi\circ T^n(x)-\phi(T^n(x)+T^n(v))|\\
&\le&\|\phi\|_{\H^\alpha_s}|T^n(v)|^\alpha\\
&\le&\|\phi\|_{\H^\alpha_s}\lambda^{-n}|v|^\alpha.
\end{eqnarray*}

In the same way, we get
\begin{equation}\label{inn2}
\|\phi\circ T^{-n}\|_{\H^\alpha_{u,c}}\le Cn^r\|\phi\|_{\H^\alpha_{u,c}}.
\end{equation}

\medskip

Now, to prove the proposition, we apply what precedes to 
$$
\phi=\prod_{j=1}^m\varphi\circ T^{k_j}\mbox{ and }\psi=\prod_{j=1}^p\varphi\circ T^{l_j}
$$
with negative $k_i$ and positive $l_i$. 
Using (\ref{inn1}) and (\ref{inn2}), the computation shows that
$$
\|\phi\|_{\H^\alpha_{u,c}}\le\sum_{j=1}^m\|\varphi\|_\infty^{m-1}\|\varphi\|_{\H^\alpha_{u,c}}|k_j|^r\le \|\varphi\|_{\H^\alpha}\sum_{j=1}^m|k_j|^r
$$
and
$$
\|\psi\|_{\H^\alpha_s}\le\sum_{j=1}^p\|\varphi\|_\infty^{p-1}\|\varphi\|_{\H^\alpha_s}\le \|\varphi\|_{\H^\alpha}
$$
Then, by (\ref{inn}), the proposition is proved.
\cqfd

Further the classical central limit theorem holds for any H\"older functions of the torus \citeaffixed{Leo60,Leb99}{see}(See Leonov \citeyear{Leo60}; Le Borgne \citeyear{Leb99}). Then Theorem~\ref{thm1} applies and we get Theorem \ref{thmtorus}.

\paragraph*{Acknowledgement}\

The authors would like to thank the referee for her/his very careful reading
of the manuscript and for several  thoughtful comments that helped improving the paper. 

\small

\bibliographystyle{agsm}

\end{document}